\documentclass[preprint,12pt]{elsarticle}

\usepackage{amssymb}
\usepackage{hyperref}
\usepackage{amsmath}
\usepackage{fixltx2e}
\usepackage{footnote}
\usepackage{multicol}
\usepackage{booktabs}
\usepackage[flushleft]{threeparttable}
\usepackage[font=small,labelfont=bf]{caption}

\newdefinition{rmk}{Remark}
\newproof{proof}{Proof}

\journal{---}

\begin{document}

\begin{frontmatter}


 \title{Algorithms for Computing Wiener Indices of Acyclic and Unicyclic Graphs}
\author{Khawaja Muhammad Fahd}
\ead{muhammad.fahd@gmail.com}
\address{Department of Computer Science, Lahore University of Management Sciences (LUMS),Lahore, Pakistan. }
\address{Department of Computer Science, Riphah Institute of Computing and Applied Sciences (RICAS),\\
	Riphah International University, 14 Ali Road, Lahore, Pakistan.\fnref{label3}}

\author{Muhammad Kamran Jamil\fnref{label2}}
\ead{m.kamran.sms@gmail.com}
\fntext[label2]{Corresponding author}
\address{Department of Mathematics, Riphah Institute of Computing and Applied Sciences (RICAS),\\
Riphah International University, 14 Ali Road, Lahore, Pakistan.\fnref{label3}}



\address{}

\begin{abstract}
	Let $G=(V(G),E(G))$ be a molecular graph, where $V(G)$ and $E(G)$ are the sets of vertices (atoms) and edges (bonds). A topological index of a molecular graph is a numerical quantity which helps to predict the chemical/physical properties of the molecules.  The Wiener, Wiener polarity and the terminal Wiener indices are the distance based topological indices. In this paper, we described a linear time algorithm {\bf(LTA)} that computes the Wiener index for acyclic graphs and extended this algorithm for unicyclic graphs. The same algorithms are modified to compute the terminal Wiener index and the Wiener polarity index. All these algorithms compute the indices in time $O(n)$.
\end{abstract}


\end{frontmatter}


\section{Introduction}

Mathematical chemistry is a branch of theoretical chemistry. This branch is used to discuss and predict the molecular structure with the help of mathematical methods without necessarily referring quantum mechanics. Recent research in mathematical chemistry, particular attention has been paid to so-called topological indices. A topological index is a number associated with the graph obtained from the chemical compound.\\

Wiener index is one of the oldest and widely used topological indices. Wiener introduced this topological index while he was working on the boiling point of alkane molecules \cite{w}. After that Wiener index is used to predict many other physical and chemical properties of the molecules such as critical constants, density, viscosity, surface tension, ultrasonic sound velocity \cite{bmpt,pdb,rc,rtn,st}. \\

Let $G=(V(G),E(G))$ be a molecular graph, where $V(G)$ and $E(G)$ represent the sets of vertices (atoms) and edges (bonds) of the graph $G$. The distance between two vertices $u$ and $v$ is the length of the shortest path between them, and is denoted as $d(u,v)$. For a vertex $v\in V$, the degree of $v$ is the number of the adjacent vertices to $v$. It is denoted as $d(v)$.\\

The Wiener index of a molecular graph $G$ is the sum of distances between all pair of vertices of $G$ i.e.,
\[W(G)=\sum_{\{u,v\}\in V(G)}d(u,v)\]

The Wiener polarity index of a molecular graph $G=(V(G),E(G))$ is defined as the number of unordered pairs of vertices \{u,v\} of $G$ such that the shortest distance $d(u,v)$ between $u$ and $v$ is $3$,
\[W_P(G)=|\{\{u,v\}|d(u,v)=3,u,v\in V(G)\}|\]
The terminal Wiener index $TW(G)$ of a molecular graph is defined as the sum of distances between all pairs of its pendent vertices. If $V_T=\{v_1,v_2,\cdots,v_k\}$ is the set of all pendent vertices of $G$, then
\[TW(G)=\sum_{\{u,v\}\subseteq V_T(G)}d(v_i,v_j)\]

A numerous work has been done on these topological indices. In 1988, Mohar et. al. \cite{mp} introduced an algorithm to find the Wiener index of acyclic graphs. Azari et. al. determined the Wiener index of molecular graphs made by hexagons \cite{ait}. Deng investigate the Wiener index of Spiro and polyphenyl hexagonal chains \cite{d}. In 2007, Ashrafi et. al. \cite{ay} gave the algorithm for computing distance matrix and Wiener index of zig-zag polyhex nanotubes. Ili\`c et. al. \cite{ii} described some algorithms for computing topological indices of chemical graphs. In this paper, we described linear time algorithms to compute the Wiener index, the Wiener polarity index and the terminal Wiener index for acyclic graphs and unicyclic graphs.

\section{The Wiener indices of acyclic and unicyclic graphs}

The main results of this section are linear time algorithms for computing Wiener index (WI) for an arbitrary acyclic and unicyclic graphs. Linear time algorithm (LTA) for acyclic graphs improves and simplifies a nice recursive algorithm of Canfield et. al. \cite{crr}. LTA is similar to the algorithm of Mohar et. al. \cite{mp} in the way that it calculates the number of shortest paths by looking at the contribution of each edge, and summing all these values. LTA achieves this without converting the tree into a rooted tree. We extend LTA to compute the Wiener Index of unicyclic graphs.

\subsection{\bf{Algorithm for Wiener index}}
\vspace{.5cm}
Following variables are used in our LTA.\\

\noindent
WI is the Wiener Index\\
L[i] is the adjacency list of vertex $i$.\\
d[i] is the degree of vertex $i$.\\
Stack is the stack of leaf nodes.\\
Is-empty() is a function for stack which returns true if stack is empty and false otherwise. Pop() is a function of stack and it returns the element that was pushed last onto the stack. Push(i) is a function of stack that pushes vertex $i$ onto the stack.\\
A vertex may be a leaf vertex in the original graph or it becomes a leaf vertex after successive deletion of leaf vertices. Degree of a leaf vertex is 1. Also every edge in a tree is a cut edge. Consider a vertex $v$ at the point it becomes a leaf vertex and its edge $e$. Removing $e$ from the original graphs disconnects the graph {\it count [v]} will give the number of vertices in the component containing $v$ and there will be a total of {\it count $[v]\times$ (n-count [v])} paths (shortest paths) going through the edge $e$ in the original graph.\\

LTA for acyclic graphs have the following steps.\\

\noindent
$V = \{1, 2, ..., n\}$\\
WI = 0\\
FOR (i = 1 to n)\\
\hspace*{12 mm} count[i] = 1\\
\hspace*{12 mm} traverse L[i] and count the number of neighbors of i and update d[i]\\
\hspace*{12 mm} IF $d[i] = 1$\\
\hspace*{25 mm} stack.push(i)\\
WHILE (stack is not empty)\\
\hspace*{12 mm} leaf = stack.pop()\\
\hspace*{12 mm} $WI = WI + count[leaf] * (n-count[leaf])$\\
\hspace*{12 mm} $d[leaf] = 0$\\
\hspace*{12 mm} traverse L[leaf] and look at its neighbor nb.\\
\hspace*{12 mm} IF $d[nb] > 0$\\
\hspace*{25 mm} $count[nb] = count[nb] + count[leaf]$\\
\hspace*{25 mm} $d[nb] = d[nb] - 1$\\
\hspace*{25 mm} IF $d[nb] = 1$\\
\hspace*{37 mm} $stack.push(nb)$\\

At the end of this algorithm, WI gives the value of the Wiener Index of given tree.\\

For unicyclic graphs, the following lines of code should be added after the algorithm described above.\\

\noindent
$i = 1$\\
while($d(i) = 0$)\\
\hspace*{12 mm} $i = i + 1$\\
$k = 0$\\
$start = i$\\
$current = start$\\
$prev = -1$\\
DO\\
\hspace*{12 mm} $val[k] = count[current]$\\
\hspace*{12 mm} $k = k + 1$\\
\hspace*{12 mm} traverse L[current] until its neighbor having $d[nb] = 2$ and $nb \neq prev$ is found\\
\hspace*{12 mm} $prev = current$\\
\hspace*{12 mm} $current = nb$\\
WHILE ($current \neq start$)\\
$diff[1] = 0$\\
FOR ($i = 1$ to $\lfloor k/2 \rfloor$)\\
\hspace*{12 mm} $diff[1] = diff[1] + val[i]$\\
FOR ($i = 2$ to $k - 1$)\\
\hspace*{12 mm} $diff[i] = diff[i-1] - val[i-1] + val[(i+\lfloor k/2 \rfloor)$ mod k]\\
$add = 0$\\
FOR ($i = 1$ to $\lfloor k/2 \rfloor$)\\
\hspace*{12 mm} $add = i * val[i]$\\
$WI = WI + val[0] * add$\\
FOR ($i = 2$ to $k - 1$)\\
\hspace*{12 mm} $add = add - diff[i] + (\lfloor k/2 \rfloor* val[(i+\lfloor k/2 \rfloor)$ mod k]\\
\hspace*{12 mm} $WI = WI + val[i] * add$]\\
subtract = 0\\
IF (k is EVEN)\\
\hspace*{12 mm} FOR( $i = 0$ to $k/2 - 1$)\\
\hspace*{25 mm} $subtract = subtract + val[i] * val[i + k/2$]\\
\hspace*{12 mm} $WI = WI - subtract$\\

Here WI gives the Wiener Index for unicyclic graphs.\\

We compared $LTA$ with Slow-All-Pairs (SAP), Faster-All-Pairs(FAP), Floyd-Warshall(FW), Breadth First Search (BFS) algorithms by taking different graphs with various number of vertices and took their average. The details are illustrated in Table \ref{Table}.

\begin{table}[h!]
	\begin{tabular}{|c|c|c|c|c|c|}
	\hline
	n&SAP  &FAP  &FW  &BFS  &LTA  \\
	\hline
	10& 0.108252 & 0.105254	  & 0.0942775	  &0.093166	     &0.084932	       \\
	\hline
	20&	0.11424&	0.109358&	0.107474&	0.10362&	0.086898\\
	\hline
	50&	0.13496&	0.13132&	0.128242&	0.116254&	0.088754\\	
	\hline
	100&	0.71328&	0.1891&	0.16542&	0.12597&	0.09881\\
	\hline
	1000&	6158&	76.15&	7.7228&	0.13312&	0.105416\\
	\hline
	10000& & & &				2.7072&	0.11802\\
	\hline
	\end{tabular}
	\caption{Comparison of five algorithms. The times are given in miliseconds and are averaged over various graphs. }\label{Table}
\end{table}

\subsection{\bf{Analysis of Our Algorithm}}
	When we delete a leaf vertex, $v\in V(G)$ its corresponding edge $uv \in E(G)$ is also deleted. At this point the number of shortest path going through the edge $uv$ is added towards our final solution but all the paths from vertex $v$ is not added. To accommodate for the missing values we are using the concept of count and the count of $v$ is added to count of $u$. Here count of $u$ represents the number of vertices including $u$ whose total distances depends upon the distances of $u$ is with the remaining vertices. We need to do number of steps proportional to the number of edges, hence LTA runs in linear time.\\
	
	We extend this idea to unicyclic graphs. We continue removing leaf vertices until we are left with a cycle. Suppose that the cycle is of order $k$, and $count[v_i]=c_i,i=1,2,\cdots,k-1$ the value of these counts are $c_0,c_1,\cdots, c_{k-1}$. We look at the distances clockwise in order to avoid duplication. We need to calculate the following term
	\[\sum_{i=0}^{k-1}\sum_{j=(i+1)mod \,k}^{(i+\lfloor\frac{k}{2}\rfloor)mod\,k}c_ic_j\big[(j-i)mod\,k\big]\]
	
	We first calculate $\sum_{j=1}^{\lfloor\frac{k}{2}\rfloor}jc_j$ in O(k) steps then we calculate the remaining $\sum_{j=i+1}^{i+\lfloor\frac{k}{2}\rfloor}c_j(j-i)$ terms in O(k) by calculating each of these $k-1$ terms in $O(1)$ steps. This is done by using the fact that
	\begin{align*}
		&\sum_{j=i+2}^{(i+\lfloor\frac{k}{2}\rfloor+1)mod\,k}c_j\big[(j-i)mod\,k\big]-\sum_{j=i+1}^{(i+\lfloor\frac{k}{2}\rfloor)mod\,k}c_j\big[(j-i)mod\,k\big]\\
		&=\lfloor\frac{k}{2}\rfloor c_{(i+\lfloor\frac{k}{2}\rfloor+1)mod\,k}-\sum_{j=i+1}^{(i+\lfloor\frac{k}{2})mod\,k}c_j
	\end{align*}
	
	In the end if $k$ is even then we have calculated the values
	\[\sum_{i=0}^{\frac{k}{2}-1}\frac{k}{2}c_ic_{i+\frac{k}{2}}\]
	multiple times. Therefore we need to fix this double counting by subtracting this value.

\subsection{ \bf{Algorithm for Terminal Wiener Index}}	
	
	Now we extend these algorithms to calculate the Terminal Wiener Index.\\
	
	Terminal Wiener Index is defined as the sum of distances between leaf vertices as opposed to the Wiener Index that is defined as  the sum of distances between all the vertices. Therefore the distances of non-leaf vertices should not be included in this algorithm. We modify LTA by initializing the count of non-leaf vertices to $0$. Hence the contribution of all such vertices are nullified.
	
	We need another modification to complete the algorithm for calculating Terminal Wiener Index correctly. We calculate and update the terminal Wiener index by considering $count[v]$ of the leaf vertex $v$ and adding the value $count[v]\times(n-count[v])$, where $n$ is the number of vertices in the graph. Since we want to calculate the distances between leaf vertices only and not all the vertices, $n$, we change this line to $count[v]\times(nleaf-count[v])$ where $nleaf$ is the total number of leaf nodes in the graph.

	With only these changes, our algorithm will find the Terminal Wiener Index of acyclic graphs as well as Terminal Wiener Index of unicyclic graphs. Since the number of steps are still the same as LTA, this algorithm also takes linear time.\\
	
\subsection{\bf{Algorithm for Wiener Polarity Index}}
	Wiener polarity index is defined as the number of pair of vertices at distance 3. Here we do not add the distances. We again update LTA to calculate Wiener Polarity Index.\\
	
	In LTA we repeatedly delete leaf vertices. When vertex v becomes leaf vertex, it has only only edge uv. Here $count[v]$ gives us the value of total number of vertices in the component containing $v$ if the edge $uv$ is deleted from the original graph. Accordingly we add $count[v]\times(n-count[v])$ to the value of Wiener Index. LTA does not look at the actual distance. In Wiener polarity index we are interested in counting pair of vertices at distance of $3$. Instead of the variable $count[v]$, we use $count1[v]$ and $count2[v]$. Here, when a vertex v becomes a leaf vertex, then after deleting its edge $uv$, $counti[v],i=1,2$ gives the number of vertices that are at distance $i$ from the vertex $v$ and is on the component containing $v$.\\
	
	When we delete the edge $uv$, instead of adding $count[v]\times(n-count[v])$ to the Wiener Index, we need to add the number of pairs that are at distance $3$ and update the values of $count1[u]$ and $count2[u]$ accordingly. We include the following lines instead of $WI = WI + count[v] \times (n-count[v])$.
	
	$WI = WI + count2[u] + count2[v] + (count1[u] \times count1[v])$\\
	$count1[u] = count1[u] + 1$\\
	$count2[u] = count2[u] + count1[v]$
	
	This calculates the Wiener Polarity Index of acyclic graphs in linear time as the number of steps are similar to the LTA.
	
	For calculating the Wiener polarity index of unicyclic graphs, we add the following values to WI after only a cycle is left.
	
	$WI = WI + \sum_{i=0}^{k} (count1[i] + count1[(i+2) mod k] + count2[i] + count2[(i+1) mod k] + (count1[i] \times count1[(i+1) mod k])$
	
	Clearly this can be calculated in linear time, therefore our algorithm calculates Wiener Polarity Index for acyclic graphs and unicyclic graphs in deterministic linear time.\\

\end{document}